# Optimal fast charging station locations for electric ridesharing with vehicle-charging station assignment

**Tai-Yu Ma*[1] and Simin Xie[2]**

[1] Luxembourg Institute of Socio-Economic Research (LISER), Department of Urban Development and Mobility, Esch-sur-Alzette, Luxembourg

[2] Katholieke Universiteit Leuven, Oude Markt 13, Leuven 3000, Belgium

*Corresponding author

## Abstract

Electrified shared mobility services need to handle charging infrastructure planning and manage their daily charging operations to minimize total charging operation time and cost. However, existing studies tend to address these problems separately. A new online vehicle-charging assignment model is proposed and integrated into the fast charging location problem for dynamic ridesharing services using electric vehicles. The latter is formulated as a bi-level optimization problem to minimize the fleet's daily charging operation time. A surrogate-assisted optimization approach is proposed to solve the combinatorial optimization problem efficiently. The proposed model is tested on a realistic flexible bus service in Luxembourg. The results show that the proposed online charging policy can effectively reduce the charging delays of the fleet compared to the state-of-the-art methods. With 10 additional DC fast chargers installed, charging operation time can be reduced up to 27.8% when applying the online charging policy under the test scenarios.

# 1. Introduction

Electrifying shared vehicle fleets has become increasingly popular due to the benefits of the zero-emission nature of battery electric vehicles (EVs) and lower per-mile operating costs (George and Zafar, 2018). However, operating EV fleets for shared on-demand mobility services requires recharging vehicles several times a day, mainly relying on DC fast chargers to reduce the idle time of vehicles (Jenn, 2019). Due to the limited number of fast chargers and the increasing number of EVs, EV charging operations might suffer increasing queuing delays and higher charging operation costs. Transport network companies need to consider new fast charging infrastructure investments to reduce the total system operating costs. Several studies have begun working on related issues for fast charging station planning of on-demand mobility services using EVs (Wu and Sioshansi, 2017; Jung and Chow, 2019; Roni et al., 2019).

Charging station location problems have been studied widely in recent years, focusing on the optimal configuration of public charging infrastructure to meet the charging needs of private electric vehicles (PEVs) and promote PEV adoption (see the recent literature review in Shen et al., 2019). Nevertheless, these studies are mainly conducted from a PEV charging needs perspective (overnight charging, trip chaining for daily activity realizations), which is characterized by very different driving and charging demand patterns of commercial electric vehicles for shared or on-demand mobility services. From the operator perspective, the most important obstacles for e-fleet operation are related to the long charging time, which might significantly reduce vehicle availability for serving customers (George et al., 2018). Operating EVs for shared mobility services requires optimal charging management given limited fast charging facilities. While there are several studies related to public charging station location planning for electric car-sharing systems and e-taxi fleets (Jung et al., 2014; Asamer et al., 2016; Brandstätter et al., 2017; Yang et al., 2017; Zhang et al., 2019), they do not jointly address the optimal charging management of e-fleets when planning new charging station infrastructure. This issue is particularly relevant because the location of charging stations influences the access time and charging queuing delays; in turn, the charging management strategy (determining at which charging stations to charge EVs over time) impacts the occupancy of charging stations and queuing delays when arriving at charging stations. The joint consideration of optimal charging management and new charging infrastructure planning allows for the optimization of e-fleet operation, minimizing the total queuing delays of charging operations under stochastic customer demand and driving patterns of EVs.

We propose a new dynamic vehicle-charging station assignment policy, formulated as a mixed integer linear program (MILP), which is integrated into the bi-level optimal fast charging location problem for fast charging station extension of dynamic ridesharing services. The vehicle-charging station assignment policy is considered in a stochastic environment accounting for vehicle charging delays (including travel time to reach charging stations, charging queuing times, and charging time) under a heterogeneous capacitated charging network. We apply the MILP model based on the rolling time-window framework to minimize the daily charging operation times of EV-based dynamic ridesharing services. To validate the proposed model, we test it on the Luxembourg flexible bus service under different hypothetical scenarios. The results show that the proposed method significantly reduces vehicle charging operation time and queuing delays compared to the state-of-the-art charging policies. To solve the bi-level fast charging location problem, we propose a surrogate-assisted optimization (SO) approach to efficiently find near-optimal solutions. The benefit of the fast charging station extension for fleet charging operation time and the environmental benefit of electrifying the fleet are evaluated for the Luxembourg case study.

The rest of the paper is organized as follows. Section 2 gives a short review of the literature related to charging station location problems and charging scheduling under uncertainty. Section 3 introduces the bi-level optimization-simulation framework for fast charging location planning by integrating a new vehicle-charging station assignment model to minimize the total vehicle idle time. An SO approach is proposed to solve the bi-level optimization problem efficiently. Section 4 applies the methodology for a realistic flexible bus service in Luxembourg. We describe the simulation platform and empirical data and conduct several experiments to evaluate the performance of the proposed methodology. We quantify the benefit of new fast charging stations on total vehicle charging operation times and evaluate the environmental benefits of the electrification of the flexible bus service in Luxembourg. Finally, conclusions are drawn and future extensions are discussed.

## 2. Related work

Charging station location problems have been widely studied in recent years in order to plan charging infrastructure to meet EV charging needs. The problem involves the optimization of the configuration (number, type, and locations) of a number of charging facilities where an objective function with EV access cost, queuing delays, and/or investment cost is minimized. The problem is generally modeled as a p-median problem to locate a number of facilities with minimal total customer access time (Serra and Marianov, 1998; Drezner and Hamacker, 2002). The basic p-median charging station location model has been extended by considering multiple types of chargers (Wang and Lin, 2013) and different factors of uncertainty. Several studies incorporate stochastic PEV driving patterns for DC fast charging station location planning. Davidov and Pantoš (2017) proposed a stochastic optimization model for long-term charging infrastructure expansion planning considering the stochastic driving behavior of PEV users under investment cost and quality of service constraints. Wu and Sioshansi (2017) proposed a stochastic optimization model for locating public fast charging infrastructure with the objective of maximizing the expected number of PEVs served by the charging infrastructure considering trip-chaining behavior. Liu and Wang (2016) proposed a tri-level optimization-simulation model for the multi-type charging station location problem considering a user-equilibrium model. The upper level considers a multi-type charging station location problem for social welfare maximization. The middle level considers user's EV purchase type choice behavior for plug-in or wireless charging EVs. The lower level considers tour-based travel demand and charging needs under a user-equilibrium model. Due to the high computational time needed for solving the lower-level problem, the authors apply an SO approach to solve the tri-level optimization problem. A recent review of the charging-station location model can be found in Shen et al. (2019).

For commercial e-fleet charging network planning, several charging-location optimization applications have been proposed for electric car-sharing services (Brandstätter et al., 2017), e-taxis (Asamer et al., 2016; Jung et al., 2014), and shared autonomous vehicles (Zhao et al., 2019; Zhang et al., 2019). However, these studies consider the charging station planning problem by neglecting the joint charging management problem. Jung et al. (2014) classify the EV charging station location models into three categories: node-based, flow-interception-based, and itinerary-interception-based approaches. The node-based approach minimizes the total weighted access travel distance/costs from node demand to the nearest facilities. The flow-interception approach determines the facility locations to maximize the coverage of origin-destination (OD) flows of electric vehicles. The itinerary-interception approach considers PEV drivers' trip-chaining behavior to minimize the total travel time/cost and delays of recharging when siting a set of charging stations. The authors proposed an itinerary-interception-based model for charging station location optimization of a fleet of e-taxis considering two new aspects: (1) stochastic charging demand and (2) queuing delays at charging stations due to daytime charges. The charging behavior of drivers is assumed using the nearest charging station policy for recharging when the battery level of vehicles is lower than a threshold. Such a non-coordinated charging policy might incur significant charging delays when multiple EVs head for the same fast charging station in rush hours (Yuan et al., 2019). With the advance of communication technology, the operator can deploy centralized control strategies to optimally assign vehicles to charging stations based on real-time monitoring of vehicle battery state and charging station queues (Hu et al., 2014). Hu et al. (2016) provide a recent literature review on the different control strategies available for fleet operators to optimize EV charging scheduling and management.

The long charging time of EVs and limited fast charging infrastructure may increase the charging delays, idle time, and operation costs of e-fleets. Recent studies have started looking at smart charging management under uncertainty to minimize charging delays for fleet operations (Qin and Zhang, 2011; Yuan et al., 2019; Tian et al., 2016; Ma et al., 2019a; Ma, 2020; Zhang and Chen, 2020; Pantelidis et al., 2020). Qin and Zhang (2011) introduced a charging scheduling model to minimize charging waiting times, considering a highway road network with charging stations deployed at some entrance/exit locations. The charging locations are connected by the Internet and communicate the charging station state with nearby EVs via wireless devices. The problem aims to minimize the average charging delay by considering an M/M/c queuing system with hypothetical Poisson-distributed service time and arrival rates. Yuan et al. (2019) proposed a partial charging policy for e-taxis to minimize the idle time of vehicles and increase the availability of vehicles to serve customers. They introduce a receding horizon optimization approach to scheduling e-taxi charging by sensing the state of charge of the e-taxi fleet and the occupancy of charging stations. The system is zone-based to assign e-taxis from one zone to another over pre-defined charging time decision slots to minimize the total vehicle idle time of the e-taxi fleet. The queuing delay at the charging station is approximated at the zonal level. Tian et al. (2016) consider a real-time charging station recommendation system for e-taxi operations to minimize the expected charging waiting time. The system tracks real-time e-taxi states (location, battery level, and driving patterns) to estimate an e-taxi's probabilistic charging intention and then recommends the station with the shortest charging waiting time. However, the recommendation system does not coordinate the charging scheduling of the entire e-taxi fleet.

Charging management allows the fleet operator to reduce the total charging delays and operation costs of on-demand or shared mobility services. However, existing optimal charging station location models have not integrated dynamic charging management into fast-charging station location decisions. This research gap is addressed in this study. We summarize the main contributions as follows:

1) We propose a dynamic vehicle-charging station assignment model for the electric ridesharing system to minimize total vehicle idle times for recharging under a stochastic environment. We explicitly consider charging station capacity and vehicle charging queuing delays to capture total charging operation times under stochastic demand. Our realistic case study shows that the proposed methodology outperforms state-of-the-art methods.

2) We integrate the above charging management policy into the facility location problem of fast charging station extensions for the electric ridesharing system. The problem is formulated as a bi-level optimization problem and solved by the surrogate-assisted optimization approach.

3) The proposed methodology is applied to a realistic flexible bus service in Luxembourg for which the benefit of new fast charging stations is evaluated. The environmental benefit of electrifying the fleet of flexible bus services in terms of $CO_2$ emission savings is assessed under different test scenarios.

## 3. Methodology

This section proposes a bi-level optimization-simulation framework for modeling optimal fast charging station location problems with integrated charging management for dynamic ridesharing systems using EVs. Section 3.1 proposes the bi-level optimization-simulation modeling framework. Section 3.2 presents the vehicle dispatching and routing policy for dynamic dial-a-ride service using EVs. An online vehicle-charger assignment model is proposed to minimize the total vehicle idle time for recharges. As the upper-level fast charging station location problem involves the computationally expensive simulation of a dynamic electric ridesharing system, an SO approach is proposed in Section 3.3 to find solutions.

3.1. Proposed bi-level optimization-simulation framework for optimal fast-charging station location

We consider the problem of siting a set of fast charging stations from a ridesharing/microtransit operator perspective (e.g., Via https://ridewithvia.com/ or MaaS Global https://whimapp.com/) under stochastic customer demand. The operator deploys centralized real-time charging management to assign EVs to specific charging stations for recharge to minimize the total charging delays. To limit the scope, we assume the number of fast chargers is given exogenously. Practitioners can relax this assumption by considering their budget constraints to jointly determine the number of fast chargers, their types, and their locations (Zhang et al., 2019) or by using a multi-period planning framework for charging infrastructure extension planning (Chung and Kwon, 2015). This type of problem might be encountered by fleet operators of shared mobility services (car-sharing, flexible buses, taxis) seeking to electrify their current gasoline vehicles and/or invest in fast charging facilities to improve the efficiency of daily charging operations. The fast-charging station location problem is formulated as a bi-level optimization problem with: (1) an upper-level facility location problem to minimize total daily vehicle idle time (travel time to reach a charging station, waiting time to be served at charging stations, and charging time) due to charging operations, and (2) a lower-level problem of dynamic ridesharing using a fleet of EVs with dynamic vehicle-charging station assignment management.

P1: Upper-level problem

$$\min Z(\boldsymbol{u}) = \sum_{v}(T_v^A(\boldsymbol{u}) + T_v^G(\boldsymbol{u}) + T_v^W(\boldsymbol{u})) \qquad (1)$$

subject to

$$\sum_{k \in K} u_k = U \qquad (2)$$

$$0 \leq u_k \leq u_k^+, k \in K \qquad (3)$$

$$u_k \in 0 \cup Z^+ \qquad (4)$$

The objective function (1) minimizes the total fleet idle time due to multiple charging operations executed during daily ridesharing service operations. The daily charging operation time of a vehicle includes the toal access time (travel time to reach a charging station) $T_v^A$, charging time $T_v^G$, and total waiting time at charging stations $T_v^W$. These metrics depend on the charging infrastructure configuration ($\boldsymbol{u}$) and are the outputs of the lower-level problem. Constraint (2) states that the total number of fast chargers to be installed is fixed. Constraint (3) states that the number of installed chargers $u_k$ at site k cannot exceed a pre-defined maximum number $u_k^+$. Constraint (4) ensures that the number of installed chargers at site k is a non-negative integer.

P2: Lower-level problem

The lower-level problem is a dynamic ridesharing problem with stochastic demand using a fleet of EVs. Due to stochastic customer arrival patterns, vehicle charging demand and queuing delays are evaluated by simulations. Existing studies assume that e-fleets recharge their battery at the nearest charging station whenever the vehicle's battery level is lower than a threshold (Bischoff and Maciejewski, 2014; Jung and Chow, 2019). The goal is to minimize the total travel time and queuing delays for daily charging operations. Different from the nearest charging station approaches, we adopt an optimal vehicle-charging station assignemnt strategy based on the current system state (queuing state at chargers, vehicle position and energy level) to minimize total vehicle idle time. In doing so, available charging capacity can be used more effectively compared to the myopic nearest charging station assignment strategy. We demonstrate the benefit of the optimal charging station assignment policy and compare its performance with state-of-the-art charging station assignment policies in Section 4. The bi-level interactions allow the iterative modeling of the impact of charging station location decisions on vehicle idle time due to charging operations at the upper level, while at the lower level, vehicle dispatching and routing are influenced by the charging operations, which provide the metrics to adjust decisions regarding charging station locations. The overall modeling framework and interactions of different modules is illustrated in Figure 1.

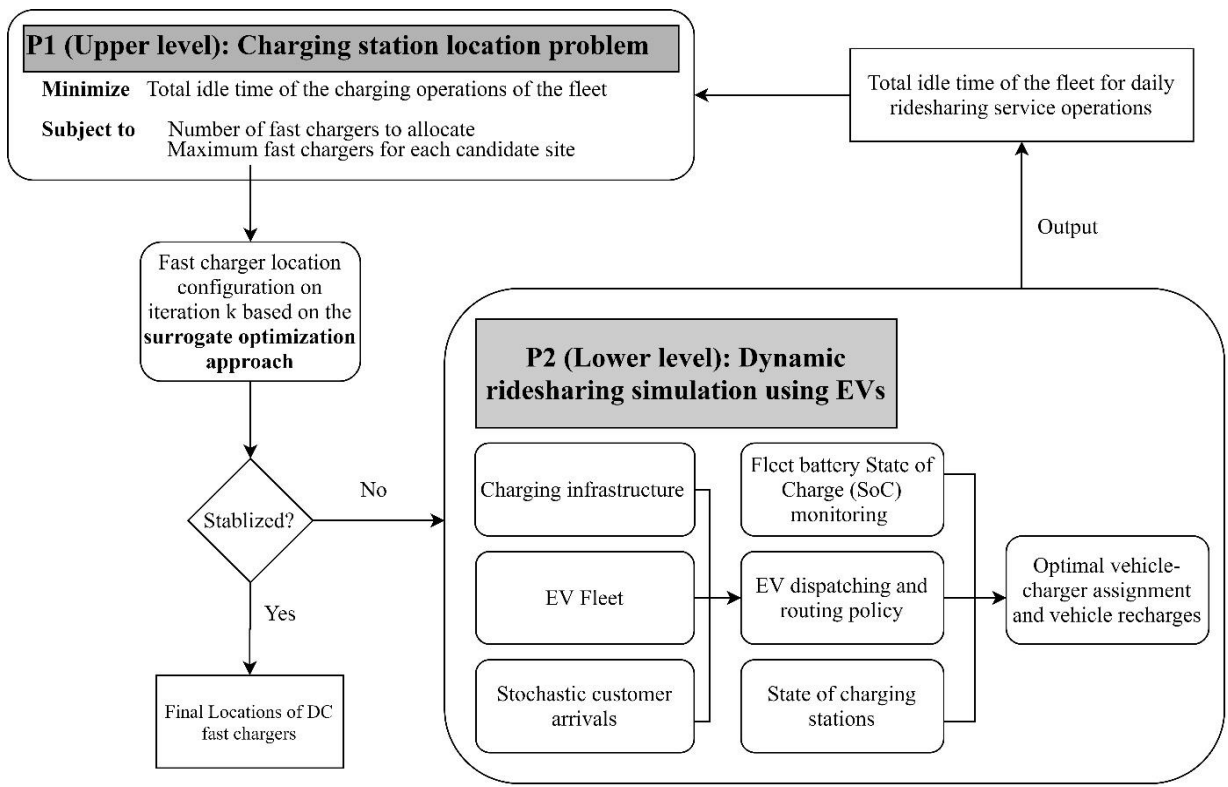

Figure 1. Overall modeling framework and interactions of different modules

3.2. Dynamic ridesharing using EVs and optimal vehicle-charger assignment

*Notation*

| | |
|---|---|
| $I$ | Set of vehicles to be assigned for recharge at the beginning of epoch $h$ (index $h$ is dropped) |
| $J$ | Set of chargers |
| $t_{ij}$ | Travel time from the location of vehicle $i$ to charger $j$ |
| $d_{ij}$ | Travel distance from the location of vehicle $i$ to charger $j$ |
| $M$ | Large positive constant |
| $e_i$ | State of charge of vehicle $i$ at the beginning of epoch $h$ (index $h$ is dropped) |
| $e_{min}$ | Minimum reserve battery level of EVs |
| $e_{max}$ | Maximum recharge level of battery |
| $t_j^A$ | Time until which charger $j$ is occupied by another vehicle from the beginning of epoch $h$ (index $h$ is dropped) |
| $\varphi_j$ | Charging power of charger $j$ (kWh) |
| $\mu$ | EV driving efficiency (kWh/km) |

*Decision variable*

| | |
|---|---|
| $X_{ij}$ | Vehicle $i$ is assigned to charger $j$ for recharge in epoch $h$ if $X_{ij} = 1$, and 0 otherwise (index $h$ is dropped) |
| $Y_{ij}$ | Amount of energy recharged for vehicle $i$ at charger $j$ |
| $W_{ij}$ | Artificial variable representing the waiting time of vehicle $i$ at charging station $j$ |

Consider a transport network company (operator) providing a ridesharing/microtransit service using a fleet of EVs. We assume that each EV is recharged to 80% of battery capacity to conserve EV battery life (Zhang and Chen, 2020). The state of charge of each EV is monitored in real-time by a dispatching center using dedicated remote-communication technology (Hu et al., 2014). For simplification, we assume a linear energy consumption of EVs, i.e., the energy consumption is proportional to travel distance regardless of other factors (Goeke, 2019). Moreover, when the state of charge of an EV is lower than a predefined threshold, a charging request is sent to the dispatching center for centralized vehicle-charging station assignment. We assume that the dispatching center has real-time information on charging station status (number of available chargers, charger types, charging schedule of each charger). To optimize charging scheduling with minimum vehicle idle time over the planning period (24 hours), we adopt a decomposition method that solves the charging management (assign EVs to specific chargers) in a rolling time-window manner. This decomposition methodology has been widely applied in different disciplines to solve complex optimization problem under undertainty (Chand et al., 2002).

The EV-based dynamic ridesharing system is based on our previous dynamic ridesharing simulation platform (Ma et al., 2019b) using the discrete event simulation approach in the environment of EVs considering a multiple queuing system. The discrete event simulation approach has been used for simulating carsharing and on-demand mobility services using shared electric vehicles (Jäger et al., 2017; Hu et al., 2019). The simulator considers the dynamic ridesharing problem and uses a non-myopic vehicle dispatching policy to minimize the additional cost when inserting a new request on existing routes of vehicles. The vehicle routing policy assumes all customers must be served at the cost of longer waiting times for some customers. Also we consider a taxi-like service with instant requests. The applied vehicle dispatching and routing policy can be substituted by other operation policies with requests booked in advance and time-window constraints associated with customer pickup or drop-off locations (Berbeglia et al., 2010). When a new request arrives, the dispatching center determines a list of candidate vehicles $V'$ that satisfy the following conditions: (1) energy feasibility, i.e., a vehicle is considered feasible if the remaining energy after serving all customers (including the new request) and returning to the depot is no less than a threshold $e_{min}$, and (2) recharge operation: a vehicle is feasible to serve customers if it is not scheduled for recharge after serving all onboard customers. The recharging detour before serving all customers onboard is not allowed. The vehicle dispatching policy is based on a non-myopic dispatching policy (Sayarshad and Chow, 2015) that assigns a new customer to the vehicle (among candidate vehicles) with the least marginal system cost as in Eq. (5).

$$\{v^*, x_t^{v^*}\} = \operatorname{argmin}_{v \in V'} [c(v, \bar{x}_t^v) - c(v, x_t^v)], \tag{5}$$

where $c(v, x)$ is a cost function with service tour $x$ of vehicle $v$ defined as Eq. (6). $\bar{x}_t^v$ is the a posteriori tour of $v$ after inserting the new request in its current tour $x_t^v$.

$$c(v, x) = \alpha T(v, x) + (1 - \alpha) \left[ \beta T(v, x)^2 + \sum_{n \in P_v} \bar{Y}_n(v, x) \right], \tag{6}$$

where $T(v, x)$ is the travel time of tour $x$. $\bar{Y}_n(v, x)$ is the total waiting time and in-vehicle travel time for all passengers $P_v$ assigned to vehicle $v$. $\alpha$ is a weight between 0 and 1 as a trade-off between operation cost and customer inconvenience. $\beta$ is a parameter between 0 and 1 to consider future approximate system delays for the current vehicle dispatching decision. The reader is referred to Ma et al. (2018; 2019b) for details regarding setting the parameter $\beta$.

For the vehicle-charger assignment, we divide the planning horizon (one day) into a set of charging decision epochs, $H = \{1, 2, \ldots, \bar{H}\}$, with the fixed time interval $\Delta$. The time interval for a decision epoch $h$ is $t_0 + h\Delta \leq t < t_0 + (h + 1)\Delta$. Vehicles reaching below the charging threshold (e.g., 20% of battery capacity) within a decision epoch $h–1$ are scheduled for recharge at the beginning of the next epoch $h$. Given the system states at the beginning of epoch $h$, the dispatching center conducts vehicle charging station assignment by solving the following MILP problem at the beginning of each decision epoch $h \in H$. The MILP problem is formulated as follows.

$$\min F(X, Y, W) = \sum_{i \in I} \sum_{j \in J} t_{ij} X_{ij} + \sum_{i \in I} \sum_{j \in J} Y_{ij}/\varphi_j + \sum_{i \in I} \sum_{j \in J} W_{ij} \tag{7}$$

subject to

$$\sum_{j \in J} X_{ij} = 1, \quad \forall i \in I \tag{8}$$

$$\sum_{i \in I} X_{ij} \leq 1, \quad \forall j \in J \tag{9}$$

$$e_{min} \leq e_i - \mu d_{ij} X_{ij} + M(1 - X_{ij}), \forall i \in I, j \in J \tag{10}$$

$$e_{max} \leq Y_{ij} + e_i - \mu d_{ij} X_{ij} + M(1 - X_{ij}), \forall i \in I, j \in J \tag{11}$$

$$Y_{ij} \leq M X_{ij}, \forall i \in I, j \in J \tag{12}$$

$$t_j^A - t_{ij} X_{ij} - M(1 - X_{ij}) \leq W_{ij}, \forall i \in I, j \in J \tag{13}$$

$$X_{ij} \in \{0,1\}, \forall i \in I, j \in J \tag{14}$$

$$Y_{ij} \geq 0, \forall i \in I, j \in J \tag{15}$$

$$W_{ij} \geq 0, \forall i \in I, j \in J \tag{16}$$

The objective function (7) minimizes total vehicle idle time due to charging operations. The first term is related to vehicle travel time to reach charging stations. The second term is the total charging times of vehicles. The third term represents vehicle waiting times at charging stations. When a charger is occupied, the next arriving vehicle needs to wait in a queue to get served. Equation (8) ensures that each vehicle is assigned to exactly one charger. Equation (9) ensures that each charger is assigned to at most one vehicle. Equation (10) ensures that the remaining energy of a vehicle should be no less than a minimum reserve value when arriving at a charging station. Equation (11) ensures that the energy level after recharge is no less than a pre-defined maximum value. Given that the objective function is to minimize total vehicle idle time, the resulting energy level after recharge is equal to $e_{max}$. Equation (12) ensures that the recharged energy amount is non-negative when a vehicle is assigned for recharge. Equation (13) states that a vehicle's waiting time at a charging station equals the difference between the time when the charger becomes available ($t_j^A$) and the arrival time of the vehicle. Note that $t_j^A$ is obtained according to the charging station occupancy state at the beginning of epoch $h$, which might not be accurate when a vehicle $i$ reaches charger $j$ at a later time $t'$ due to different uncertainty factors, i.e., vehicles previously assigned to a charger $j$ might arrive within or later than epoch $h$ (not being considered in $t_j^A$), and the arrival time of an assigned vehicle $i$ at charger $j$ might be delayed due to serving onboard customers, which might cause another vehicle to arrive at $j$ earlier and make vehicle $i$ wait in the queue. Both situations could lead to $t_j^A$ being underestimated and compromise the effectiveness of the proposed method. To overcome this issue, an alternative method is to consider only unoccupied chargers at the moment of assignment. Our numerical study shows that adopting this alternative policy allows alleviating this uncertainty issue and results in lower total charging waiting times of the fleet. Note that for simplicity, we do not model and simulate the charging demand of other EVs. Future studies can further incorporate this uncertainty into the model.

Note that Equations (8)–(9) are correct in the situation where $|I| \leq |J|$. In case the number of vehicles is greater than that of chargers ($|I| > |J|$), Equations (8)–(9) are replaced by (17)–(18) below to ensure consistency.

$$\sum_{j \in J} X_{ij} \leq 1, \quad \forall i \in I \tag{17}$$

$$\sum_{i \in I} X_{ij} = 1, \quad \forall j \in J \tag{18}$$

The MILP problem of (7)–(18) is a variant of the generalized assignment problem, which is NP-hard (Fisher and Jaikumar, 1981). The exact solutions can be obtained for small test instances (less than 50 vehicles and chargers). We apply the Lagrangian relaxation (LR) method (Fisher, 2004) to efficiently find near-optimal solutions for large-scale test instances. The applied LR heuristic (Ma , 2020) can solve large test instances up to one thousand vehicles and one thousand chargers to near-optimal (0.5% optimality gap) in less than 3 minutes on a classical laptop. The reader is referred to that study for a more detailed description. Note that we assume that EVs follow their assigned charging station instructions for recharge when all onboard customers are served. The charging demand from private EVs or taxis that compete for available public charging resources are not considered in this study. One can relax this assumption by considering historical charging patterns of the public chargers to simulate the uncertainty related to charging infrastructure availability in the presence of other EVs.

To evaluate the proposed vehicle-charger assignment policy, two widely-used reference policies are compared. The first reference policy is referred to as the **nearest charging station policy (NCP)**. This policy assumes that an idle EV goes to the nearest unoccupied charging station for recharge when its battery level is lower than a pre-defined threshold (20% of battery capacity). This policy is widely used in the state-of-the-art methods (Bischoff and Maciejewski 2014; Jung et al., 2014). The second reference policy is referred to as the **first-come-first-served (FCFS) policy**, which minimizes vehicle charging operation time based on instant current charging station state in a FCFS manner. This policy is similar to the real-time charging recommendation system for electric taxis (Tian et al., 2017), which recommends a taxi to a specific charger to minimize vehicle idle time due to recharge. Under the FCFS policy, an idle EV goes to recharge immediately whenever the battery level is lower than the predefined threshold (20% of battery capacity). Under this policy, a vehicle $v$ is assigned to charger $j^*(v,t)$ with the lowest estimated charging operation time at time $t$ as in (19).

$$j^*(v,t) = argmin_{j \in J} \left\{ t_{ij} + \frac{e_{max} - e_v(j,t')}{\varphi_j} + \widetilde{W}_{ij}(t) \right\}, \quad (19)$$

where $i$ is vehicle $v$'s current location. $e_v(j,t')$ is the battery level when arriving at charging station $j$ at time $t'$. $\widetilde{W}_{ij}(t)$ is the expected waiting time of charger $j$ based on the occupancy information of charging stations at time $t$. $\varphi_j$ is the charging power of charger $j$.

We refer to the proposed charging management policy of Eqs. (7)–(18) as the **optimal charging policy (OCP).** The variant of the OCP policy considering only unoccupied chargers at the beginning of each decision epoch is referred to as the **OCP-A policy**.

3.3. Surrogate-assisted optimization algorithm

The SO method is a derive-free global optimization method for expensive black-box objective function evaluations given a set of constraints (Regis and Shoemaker, 2007; Vu et al., 2017). Given a few initial samples, the SO method constructs a surrogate model to determine the next promising sampling points to speed up the search of a global minimum. The process is iterated until a stop criterion is met. This method has been widely applied in both academia and industry (Vu et al., 2017). Regis and Shoemaker (2007) proved that under mild assumptions, the SO method convergences to a global optimum.

The proposed SO method for solving the bi-level charging station location problem is described as follows.

**Step 1 (Initial sampling evaluation)**: Select $n_0$ initial random evaluation points $Q_0 = \{u_k\}, k = 1, \ldots, n_0$, satisfying constraints (2)–(4) of P1. Run the ridesharing system simulations using EVs, solve the optimal charging assignment for each charging decision epoch $h \in H$, and obtain the metrics of charging delays. Evaluate the objective function (1) and obtain the pairs of the evaluation points and the objective function values $(u_k, Z_k), k = 1, \ldots, n_0$. Set the previously evaluated points as $Q = Q_0$. Set iteration index $j = 1$.

**Step 2 (Surrogate model construction)**: Given the previously evaluated points $Q$, build/update the surrogate model based on a radial basis function (Gutmann, 2001) for the objective function approximation as in (20).

$$s(\boldsymbol{u}) = \sum_{i=1}^{n} \lambda_i \phi(\|\boldsymbol{u} - \boldsymbol{u}_i\|) + p(\boldsymbol{u}), \tag{20}$$

where $\lambda_i \in \mathcal{R}$ and $\|\cdot\|$ is the Euclidean norm. $\phi(\boldsymbol{u})$ is a radial basis function. $p(\boldsymbol{u})$ is a polynomial function, defined as $p(\boldsymbol{u}) = \boldsymbol{b}^T \boldsymbol{u} + \boldsymbol{a}$, where $\boldsymbol{a}$ and $\boldsymbol{b}$ are unknown parameters. Popular radius functions such as $\phi(r) = r^3$ (cubic) or $\phi(r) = e^{-\gamma r^2}$ (Gaussian) can be used (Gutmann, 2001). The unknown parameters $\boldsymbol{\lambda}, \boldsymbol{a},$ and $\boldsymbol{b}$ can be obtained by solving the system of linear equations (Vu et al., 2017).

**Step 3 (Select the next evaluation points)**: Sample $n_j$ random feasible evaluation points based on the stochastic response surface method (Regis and Shoemaker, 2007). First, randomly generate a set of feasible candidate points $\Omega_j$ and evaluate the score of each point based on the merit function defined by (21).

$$f(\boldsymbol{u}) = w_j A_j(\boldsymbol{u}) + (1 - w_j) B_j(\boldsymbol{u}), \tag{21}$$

where $A_j = \frac{s_j(\boldsymbol{u}) - s_j^{min}}{s_j^{max} - s_j^{min}}$ if $s_j^{max} \neq s_j^{min}$, and 1 otherwise. $B_j = \frac{d_j^{max} - d_j(\boldsymbol{u})}{d_j^{max} - d_j^{min}}$ if $d_j^{max} \neq d_j^{min}$, and 1 otherwise. $s_j(\boldsymbol{u})$ is the estimated surrogate value of (20) at the candidate point $\boldsymbol{u}$. $s_j^{min}$ and $s_j^{max}$ are the minimum and maximum surrogate values in $\Omega_j$, respectively. $d_j(\boldsymbol{u})$ is a distance metric that measures the minimum distance between a candidate point $\boldsymbol{u}$ and previously evaluated points $Q$, defined as $d(\boldsymbol{u}) = \min_{\boldsymbol{u}_i \in Q}(D(\boldsymbol{u}, \boldsymbol{u}_i))$, where $D(\cdot)$ is the Euclidean distance metric. We define $d_j^{max} = \max_{\boldsymbol{u}_i \in \Omega_j}\{d(\boldsymbol{u}_i)\}$ and $d_j^{min} = \min_{\boldsymbol{u}_i \in \Omega_j}\{d(\boldsymbol{u}_i)\}$. $w_j$ is the weight between 0 and 1. The next evaluation point at iteration $j$ is selected as $\boldsymbol{u}'_j = argmin_{\boldsymbol{u}_i \in \Omega_j}\{f(\boldsymbol{u}_i)\}$. Given $\boldsymbol{u}'_j$, run the simulation of the lower-level problem and evaluate the objective function (1). Update $Q = Q \cup \boldsymbol{u}'_j$.

**Step 4 (Stop criteria)**: If iteration $j = j^{max}$ or the value of the objective function stabilizes, then stop. Otherwise set $j := j + 1$ and go to Step 2.

## 4. Computational study

In this section, we present a realistic case study for a flexible bus service using EVs in Luxembourg. First, we present the test scenarios and parameter settings. Then we describe the vehicle dispatching and routing policy of the dynamic EV-based ridesharing service and the simulation platform. To validate the proposed optimal vehicle-charging station assignment policy, we compare its performance with two state-of-the-art methods. Finally, we evaluate the environmental benefit and charging delay reduction due to new fast charging infrastructure. The simulation case study is implemented in MATLAB using a Dell Latitude E5470 laptop with win64 OS, Intel i5-6300U CPU, 2 Cores, and 8GB memory.

4.1. Luxembourg flexible electric bus case study

We test the proposed bi-level fast charging station location extension problem on a flexible bus service (Flexibus, https://www.sales-lentz.lu/fr/communes/flexibus/) in Luxembourg. Sales-Lentz (transport network company) is the major bus company in Luxembourg, operating a fleet of gasoline shuttles to provide flexible door-to-door bus services in all of Luxembourg. To promote e-mobility, the company operates a fleet of hybrid electric vehicles and fully electric minibuses in Luxembourg. In our case study, we assume that the company would like to electrify the fleet of its Flexibus service and plans to install a number of DC fast chargers in Luxembourg to reduce total charging operation time. We assume that the number of DC fast chargers to be installed is exogenously given. As public transport in Luxembourg is completely subsidized, the objective of the operator is to minimize the total vehicle-charging operation time for its daily operation of the Flexibus service.

We assume that the electric shuttles can either be recharged on the existing public charging network (Chargy network, https://chargy.lu/en/) or new fast chargers can be installed. Presently, there are a total of 814 level-2 chargers (22 kW charging power) located at 302 locations in Luxembourg, shown in Figure 2. We assume that the fast chargers are all 50 kW DC fast chargers and that the number of fast chargers to be installed is 10. The candidate locations are on the existing Chargy network (302 locations). As we have no empirical ride data from the operator, we generate random customer demand from the recent mobility behavior survey (Luxmobil, 2017) in Luxembourg, based on the probability of trip occurrences. The Luxmobil survey, which was conducted in 2017, surveyed 40,000 households in Luxembourg and 45,000 cross-border workers, with response rates ranging from 26% to 30%. The survey data contains anonymized trip chain information and individuals' socio-demographic characteristics. We assume that the fleet size of the Flexibus service is 50 electric shuttles and that its daily customer demand on weekday is 1000, corresponding to 20 customers/vehicle/day. The practitioner is encouraged to consider one week demand patterns to take into account the variability of trip demand during weekdays and weekend. The spatial distribution of the trip requests is shown in Figure 3. Under this setting, the average customer waiting time and journey time are around 15 minutes and 36 minutes using gasoline shuttles, respectively. Note that the practitioner is encouraged to use one-week historical ride data to consider the spatial and temporal variability of demand. The distribution of customer arrivals is shown on the right side of Figure 4. We observe that a morning and afternoon peak hours are between 7:00-8:00 and 17:00-18:0, respectively. The trip length distribution is shown on the left side of Figure 4. The length of 90% trips is between 5-20km with the mean and variance of 11.9 km and 23.5, respectively.

The characteristics of electric shuttles are hypothetically based on Volkswagen's 8-seat 100% electric Tribus (https://www.tribus-group.com/e-mobility/). The battery capacity of the Tribus is 35.8 kWh, with a practical driving range up to 150 km. We assume that the energy efficiency is constant and is defined as the ratio of battery capacity and the practical driving range. The fleet is assumed to be distributed to 13 different operator depots with 4 vehicles/depot, except the municipalities of Contern and Steinfort (3 vehicles/depot), as shown in Figure 2. We assume that the depots are located in the following 13 municipalities: Bettembourg, Dudelange, Esch-sur-Alzette, Mersch, Reckange-sur-Mess, Roeser, Rumelange, Walferdange, Sanem, Garnich, Koerich, Steinfort, and Contern, based on Flexibus's current service areas. The details of the simulation parameter settings are reported in Table 1.

Table 1. The parameter settings of the Flexibus case study

| Simulation settings | | EV and charging facility settings | |
|---|---|---|---|
| Number of customers | 1000 | Energy efficiency ($\mu$)[2] | 0.2387 (kWh/km) |
| Number of DC fast chargers to be installed | 10 | Battery capacity (B)[1] | 35.8 kWh |
| Fleet size | 50 | Charging power ($\varphi_{L2}$) | 22/60 (kW/min.) |
| Capacity of vehicles[1] | 8 pers./veh. | Charging power ($\varphi_{DC\ fast}$) | 50/60 (kW/min.) |
| Vehicle speed | 50 km/hour | Number of existing level-2 public chargers | 814 |
| Planning horizon ($T$) | 6:30–22:00 | Number of depots | 13 |
| Charging decision epoch ($\Delta$)[4] | 30 min. | $e_{max}$ | 0.8B |
| $\beta$ | 0.025 | $e_{min}$ | 0.1B |
| $\gamma$ | 0.5 | Predefined threshold to go recharge | 0.2B |

Remark: 1. EV characteristics are based on a Volkswagen 100% electric-powered minibus with a 150 km range (https://www.tribus-group.com/e-mobility/). 2. Based on the driving efficiency of the Tribus ($\mu = 35.8/150$). 4. corresponding to the charging time from $e_{min}$ to $e_{max}$ of a Tribus at a DC fast charger.

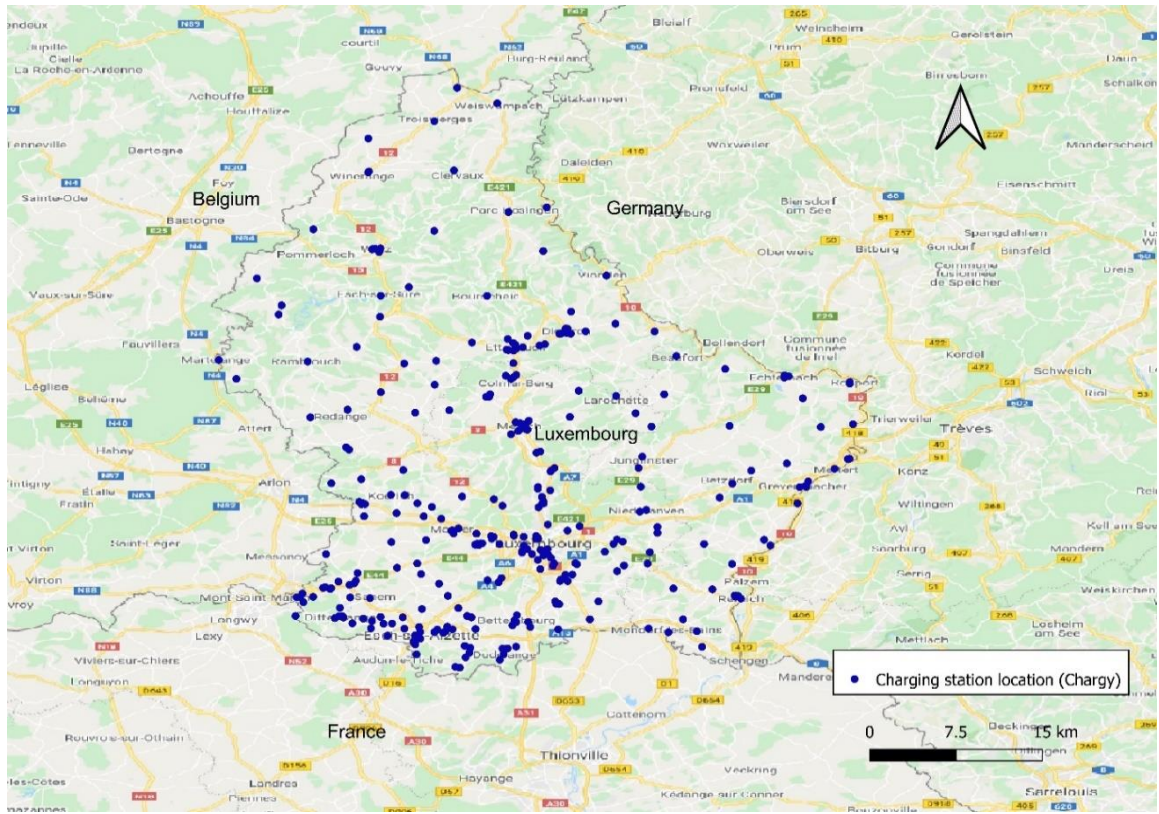

Figure 2. The spatial distribution of public charging stations in Luxembourg (Chargy network).

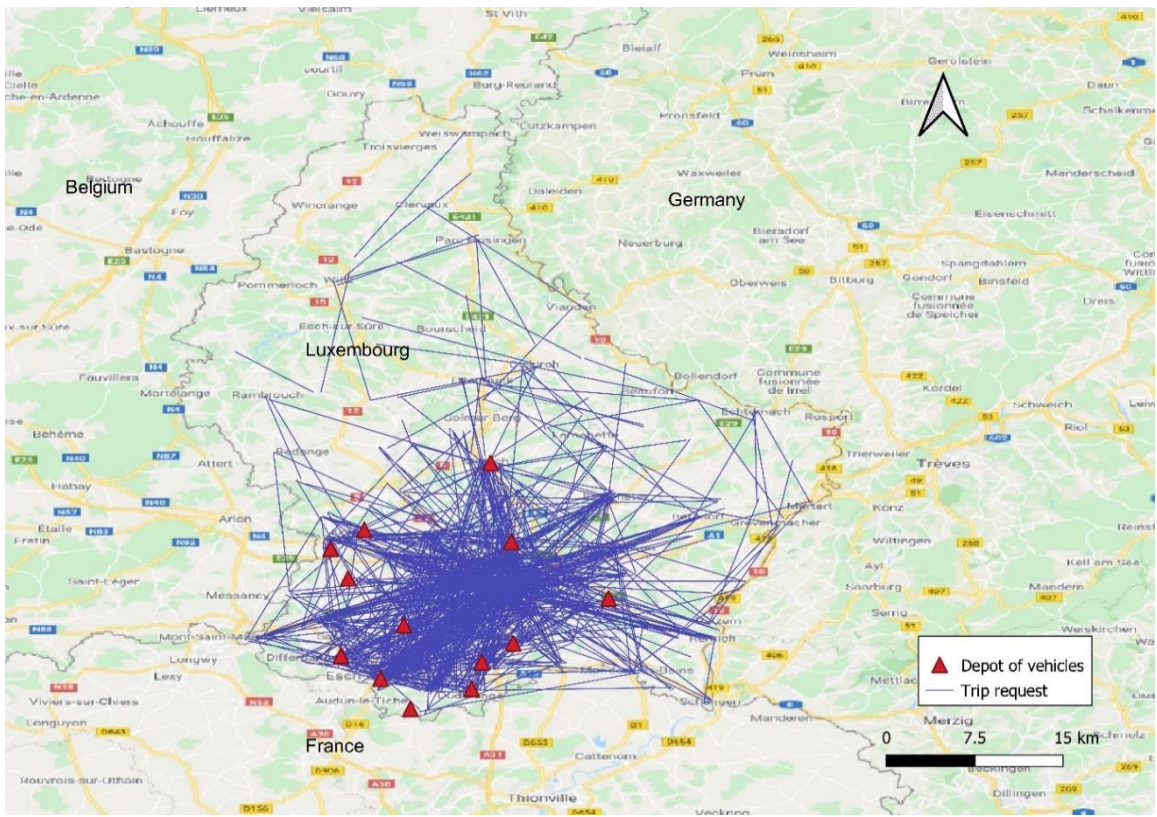

Figure 3. The spatial distribution of customer requests and vehicle depots.

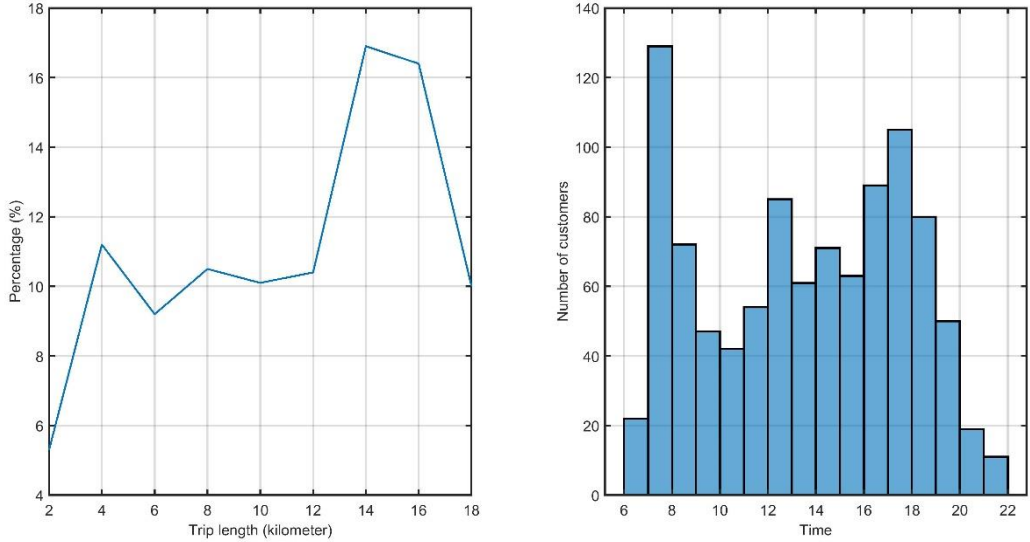

Figure 4. Example of the distribution of the length of the generated trips for the scenario of 1000 trip requests (Left); Arrival time of customers (Right).

4.2. Results

We first present the results of the proposed charging management policy to validate the proposed methodology. Then we report the computational results of the bi-level fast-charging station location problem for the Flexibus case study. Finally, we evaluate the environmental benefit in terms of annual $CO_2$ emission reduction due to fleet electrification.

4.2.1. Numerical study of the charging management policy

To take into account demand variability, we generate two data sets with 1000 and 2000 customers over 24 hours, respectively. Each dataset contains 5 test instances, generated randomly from the Luxmobil survey data based on the probability of trip occurrence. The mean and variance of trip lengths for the scenario of 2000 customers are 12 km and 23 km, respectively. The considered charging infrastructure includes the public Chargy network and 10 DC fast chargers. The locations of the 10 DC fast chargers are assumed to be around the drop-off locations of customers (Asamer et al., 2016) based on the k-means clustering method, as shown in Figure 6. We assume few charging needs from other private and commercial EVs, and thus these are ignored. Note that the optimal locations of DC fast chargers will be determined by solving the bi-level model in the next section. The results of the three charging management policies are compared in Table 2. The reported values are based on an average of 5 test instances for each demand intensity.

We observe that the OCP and OCP-A policies have lower average charging operation times per vehicle per charging operation (54–57 minutes, column (1)+(2) in Table 2) compared to the NCP and FCFS policies (65–67 minutes) with 1000 customers. The total charging waiting times of the OCP and OCP-A policies are lower compared to the NCP and FCFS policies (0–3.7 hours vs. 6–9.3 hours for one-day operations). When demand is doubled, the OCP-A policy performs best, with the lowest average charging operation time (62 minutes) compared to the other three policies (around 67 minutes for the NCP and FCFS; 76 minutes for the OCP). This shows that when charging demand is too high, the performance of the OCP policy is negatively influenced due to the inaccurate estimation of charger occupancy, resulting in additional queuing delays. We found that the OCP-A policy has almost zero charging waiting times for the fleet (0.2 hours) compared to those of the NCP (4 hours) and the FCFS (20 hours). When more fast chargers become available, the benefit of adopting the OCP-A policy for reducing charging operation times could be further increased.

In terms of customer inconvenience, we found that the OCP and OCP-A policies result in a higher mean passenger waiting time (+ around 2 minutes) and journey time (+ around 4 minutes) compared to the NCP and FCFS policies, given demand from 1000 customers. The rate of served customers is about 5% lower compared to the NCP and FCFS policies. This might be explained by some spatial mismatch effect where the k-means clustering method locates the DC fast chargers in the southern region and EVs tend to use the DC fast chargers under the OCP-based policies. Consequently, trips distributed in the far eastern and northern regions tend to be unserved (see Figure 2). When customer demand is doubled, the OCP and OCP-A policies have similar passenger waiting times and journey times compared to the other policies. The rates of served customers of the FCFS policy and the OCP-based policies are much higher (66.4%-70.6%) compared to when using the NCP policy (55.5%). The latter is due to the fact that EVs tend to use nearby level-2 chargers to recharge, resulting in higher charging times and fewer vehicles available to serve customers. We can conclude that adopting the OCP-A policy can significantly reduce charging waiting time delays and vehicle charging times in a heterogeneous charging network for both demand scenarios. In terms of the impact of the proposed OCP and OCP-A policies on customer inconvenience, this might be influenced by several factors such as the spatial distribution of customer demand, fleet configuration, vehicle driving range, and charging infrastructure configurations. Thus, we adopt the OCP-A policy as our charging management policy to determine the optimal fast charging station locations in the next section.

Table 2. System performance of different charging policies.

| Number of customers | Charging policy[1] | Charging operation | | | | Customer inconvenience | | |
|---|---|---|---|---|---|---|---|---|
| | | (1) Average charging waiting time[2] (SD) | (2) Average charging time[2] (SD) | (1)+(2) (SD) | Total charging waiting time of the fleet (hours) | MWT[3] | MJT[4] | Rate of served customers |
| 1000 | NCP | 2.4(13.0) | 65.1(2.6) | 67.5(13.2) | 6.0 | 14.5 | 35.6 | 98.6% |
| | FCFS | 3.7(14.8) | 61.6(11.3) | 65.3(15.9) | 9.3 | 14.7 | 35.6 | 98.4% |
| | OCP | 1.5(10.1) | 52.5(18.4) | 54.1(21.8) | 3.7 | 17.2 | 39.7 | 94.0% |
| | OCP-A | 0(0) | 57.1(15.8) | 57.1(15.8) | 0.0 | 16.6 | 39.0 | 93.0% |
| 2000 | NCP | 1.9(12.5) | 66.7(1.8) | 68.6(12.7) | 4.0 | 30.1 | 54.5 | 55.5% |
| | FCFS | 6.3(21) | 62.4(12.1) | 68.8(19.6) | 20.0 | 30.7 | 55.3 | 70.6% |
| | OCP | 16.7(59.9) | 59.2(15.2) | 75.9(64.0) | 50.3 | 29.8 | 55.0 | 66.4% |
| | OCP-A | 0.1(1.7) | 62.0(12.2) | 62.1(12.4) | 0.2 | 29.3 | 54.4 | 68.1% |

Remark: 1. NCP: nearest charging station policy; FCFS: first-come-first-served policy; OCP: optimal charging policy; OCP-A: optimal charging policy using unoccupied chargers only. 2. Average charging waiting time and average charging time are measured as minutes/recharge/vehicle. 3. MWT: Mean passenger waiting time (minutes). 4. MJT: Mean passenger journey time (minutes). 4 The average time of a simulation run is around 14.1 and 19.8 minutes for 1000 and 2000 customers under the OCP-A policy, respectively. 5. SD: standard deviation.

4.2.2. Optimal fast-charging station locations and the environmental benefits of the electrification of the Flexibus service

The bi-level fast-charging station location problem is solved by the SO method with the OCP-A charging management policy. We run the SO method several times with random initial solutions until the obtained solutions cannot be improved. The computational time of a simulation run is 14.1 minutes. The overall computational time for solving the optimal fast charging station location problem using the SO approach is 32.3 hours (150 simulation runs). The convergence result of the SO method is shown in Figure 5. We observe that the SO method can effectively explore the global minimum with few expensive function evaluations (lower-level simulation runs). We retain the optimal value of the objective function as 7026.63 minutes. This means that the total charging operation time for a fleet of 50 electric 8-seat minibuses to serve 1000 customers per day requires around 117 hours using exclusively the current public Chargy network and 10 DC fast chargers. However, using the k-means clustering approach for locating the DC fast chargers would result in a higher total idle time (7797.43 minutes), increasing 10.97% compared with the SO method. The spatial location of the DC fast chargers using both methods is shown in Figure 6. We observe that the SO approach locates the fast chargers around Luxembourg City, Esch-sur-Alzette, and Bettembourg, where customer demand is higher. Whereas using the k-means clustering method results in more dispersed locations of the fast charging stations on the southern half of Luxembourg."

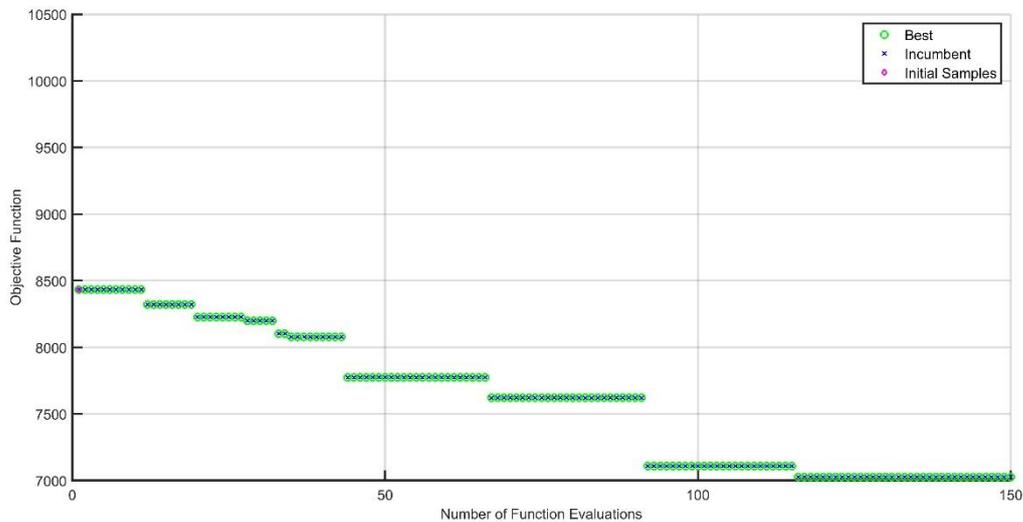

Figure 5. The convergence result of the bi-level fast-charging station location problem.

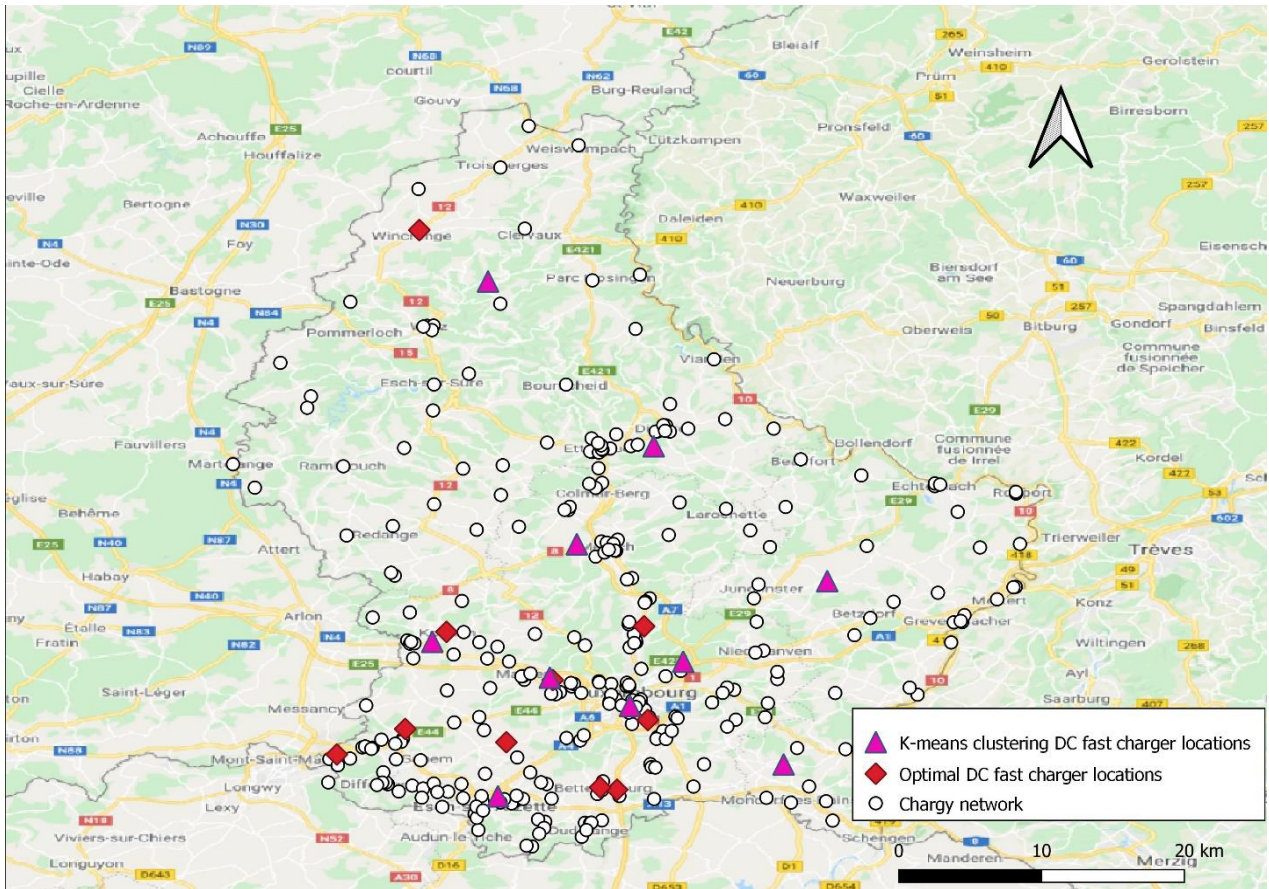

Figure 6. Spatial distributions of the optimal DC fast charger locations and that using the k-means clustering method.

We further evaluate the benefit of installing the 10 DC fast chargers in reducing total daily charging operation time (including travel time to reach chargers, waiting time, and charging time) by considering customer demand uncertainty. We vary customer demand from 500 to 2000, for which we generate 5 test instances from the 2017 Luxmobil survey data for each scenario. The result is shown on the left side of Figure 7. We found that installing an additional 10 DC fast chargers on the existing Chargy network can effectively reduce total charging operation time from 10.8% to 27.8% compared to that using the public Chargy network only. The practitioner can further conduct a sensitivity analysis to analyze the trade-off between the number of chargers and charging-time savings to plan their charging infrastructure.

In terms of the potential environmental benefit of electrifying the flexible bus fleet, we compare the annual $CO_2$ emission savings if the operator switches from its current gasoline shuttles to electric ones. We assume that operating EVs produces zero $CO_2$ emissions (Lokhandwala and Cai, 2020) while operating gasoline shuttles produces 147 grams of $CO_2$ per kilometer travelled. This estimation is based on the EU fleet-wide average emission target from 2020 onward.[1] We conduct 5 simulation runs on the test instances for each demand scenario and average the performance metrics. When daily ride requests vary from 500 to 2000, the average kilometers travelled per vehicle vary from 182 to 534 km, given a fleet size of 50 gasoline shuttles. Considering the current 6-days per week operation and 52 weeks/year, we can obtain the annual $CO_2$ emission savings when electrifying the entire fleet. The right side of Figure 7 shows that such electrification would enable the operator to reduce annual $CO_2$ emissions by 417 tons (500 customers/day) to 1225 tons (2000 customers/day). To respond to the question of how the optimal fast charger placement could reduce the $CO_2$ emission compared with systems without optimal fast charger placement, we compare the distance traveled by EVs to serve the same customer demand under systems with and without fast chargers. The result shows that installing 10 DC fast chargers allows vehicles to reduce charging operation time and serve more customers at the cost of longer vehicle travel distance per day. For the best case in terms of $CO_2$ emission savings (scenario of 2000 customers/day), an additional annual savings of 13.1 tons can be achieved based on 1389 customers served using the above computational method.

---

[1] https://ec.europa.eu/clima/policies/transport/vehicles/vans_en

Note that the above $CO_2$ savings for electrifying the flexible bus fleet do not consider the emissions generated in the electricity production phase. Depending on different types of energy generation systems, the carbon footprints vary considerably. For example, coal-powered electricity generation would generate around 500g $CO_2$eq/kWh while the carbon footprint for nuclear-powered electricity generation is around 5g $CO_2$eq/kWh (Parliamentary Office of Science and Technology, 2006). By converting the amount of the charged energy of 1781 kWh/day (500 customer scenario) and 4463 kWh/day (2000 customer scenario), 277–696 tons/year $CO_2$ emission could be generated for coal-powered energy generation and 2.77–6.96 tons for nuclear-powered energy generation.

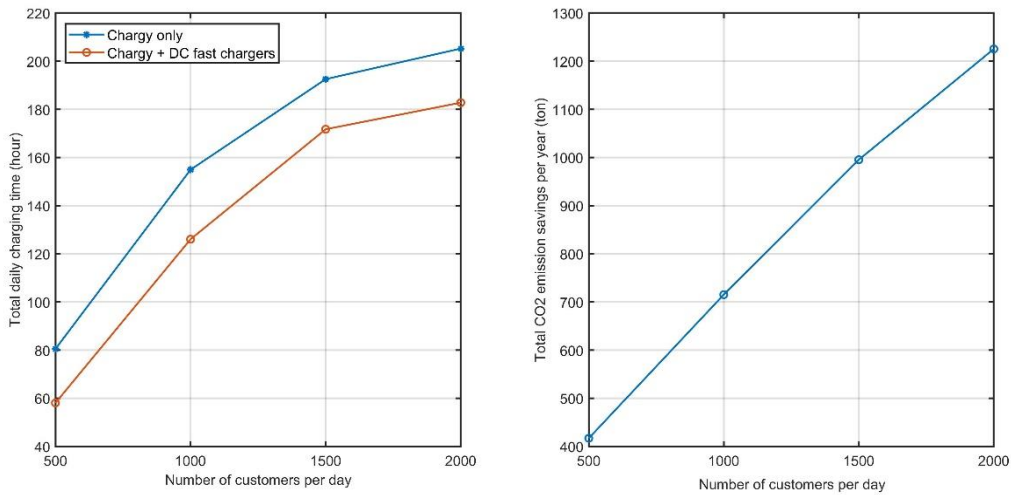

Figure 7. Environmental impact of the electrification of the flexible bus service (on the right); impact of the installation of 10 DC fast chargers on total daily charging time under different customer demands (on the left).

## 5. Discussions and conclusion

Electrifying the conventional gasoline fleet of a ridesharing service requires considering the management of daily charging operations and charging infrastructure planning to reduce charging operation delays. Existing charging station location planning methods focus mainly on private EVs or consider simple nearest charging station policies for recharge. However, none of the current studies have addressed the charging location planning problem by jointly considering optimal charging management from the fleet operator perspective.

In this study, we propose a new vehicle-charging station assignment policy to minimize total vehicle idle time due to charging operations for dynamic ridesharing services. The proposed online charging policy is integrated into a bi-level optimization-simulation model to optimize a number of fast-charging station locations under stochastic customer demand. Vehicle waiting time at charging stations is considered as a multi-server queuing system with a vehicle's stochastic charging demand and the heterogeneous characteristics of charging infrastructure. The proposed model allows optimizing charging station planning for transport network companies by integrating an optimal charging management strategy to minimize total charging times of daily operations. As the bi-level optimization problem is a difficult combinatorial optimization problem, we propose a surrogate-assisted optimization method to find good solutions.

We apply the proposed method to a realistic flexible bus service case study in Luxembourg. We find that the proposed dynamic vehicle-charging assignment policy can effectively reduce vehicle idle time and waiting times compared to the state-of-the-art charging strategies. We find that by installing 10 additional DC fast chargers at optimal locations, total charging time savings of 10.8% to 27.8% (depending on demand intensity) for a fleet of 50 (8-seater) minibuses can be obtained compared to using Luxembourg's public Chargy network only. We further evaluate the environmental benefit of electrifying the gasoline fleet, finding that the operator can be expected to reduce $CO_2$ emissions by up to 1225 tons per year given the test scenario with 2000 customers/day.

We discuss the limitation and the scalability of the proposed approach.
a) Our bi-level charging station location model minimizes the overall idle time of vehicles for daily charging operations, given exogenous customer demand, fleet size, and the number of chargers to allocate. The considered objective function can be further tailed to minimize customer inconveniences or maximize the operator's profit. The proposed model can be also extended to consider the joint optimization of the fleet size and charging infrastructure configuration (number of fast chargers, different charging powers, and charging station locations) while satisfying a limited budget constraint and customer inconvenience (Wang and Lin, 2013; Zhang et al., 2019).
b) Currently, our vehicle dispatching and routing policy does not consider time-window constraints to pick up and drop off customers. A more realistic dynamic dial-a-ride model can be designed to take into account maximum waiting times of customers, on-duty time of drivers, and maximum allowed ride/detour time for each customer (Parragh et al., 2008).
c) As for the scalability of the model, one can extend the problem in a multi-period charging station planning problem setting by incorporating a customer demand forecasting model into the annual charging infrastructure planning decision problem (Chung and Kwon, 2015). In doing so, the operator can make better decisions from a long-term perspective to anticipate future customer demand for better resource allocations.

Several research perspectives can be addressed in the future. First, the fast charging location problem can be extended to incorporate heterogeneous charger types given budget constraints. Second, one can incorporate a predictive model based on the empirical occupancy of chargers to forecast vehicle waiting times when arriving at charging stations. In doing so, the obtained charging delays reflect more realistically charger availability with the presence of other private and commercial EVs. Third, the practitioner can apply the model to planning charging infrastructure based on historical ride data over a longer period to account for demand variability. Moreover, deepened sensitivity analyses related to the impact of charging powers and battery range of EVs (Hu et al., 2019) can be conducted to evaluate the trade-offs of charging operation cost savings and charging infrastructure investment costs.

**Acknowledgements**

We thank the Ministry of Mobility and Public Works of Luxembourg for providing the anonymized 2017 Luxmobil survey data.